\documentclass{amsart}
\usepackage{latexsym, amsmath, amssymb}
\newtheorem{theorem}{Theorem}
\newtheorem{lemma}{Lemma}

\newcommand{\PA}{\par \hskip0.4cm}
\newcommand{\Z}{{\mathbb Z}}
\newcommand{\R}{{\mathbb R}}
\newcommand{\N}{{\mathbb N}}
\newcommand{\PP}{{\mathbb P}}
\newcommand{\z}{{\mathcal Z}}
\title{On the occupation measure of super-Brownian motion}

\begin{document}
\maketitle
\begin{center}
{\sc Jean-Fran\c cois LE GALL and Mathieu MERLE}

\medskip
{\small DMA-ENS, 45 rue d'Ulm, 75005 PARIS, FRANCE}
\end{center}

\begin{abstract} We derive the asymptotic behavior of the occupation measure
$\z(B_1)$ of the unit ball for super-Brownian motion started from the 
Dirac measure at a distant point $x$ and conditioned to hit the unit ball.
In the critical dimension $d=4$, we obtain a limiting exponential
distribution for the ratio $\z(B_1)/\log|x|$. 
\end{abstract}

\section{Introduction}

\PA The results of the present work are motivated by
the following simple problem about branching random walk in
$\Z^d$. Consider a population of branching particles in $\Z^d$,
such that individuals move independently in discrete time
according to a random walk with zero mean and finite
second moments, and at each integer time individuals die and give
rise independently to a random number of offspring according to
a critical offspring distribution. Suppose that the population starts
with a single individual sitting at a point $x\in\Z^d$
located far away from the origin, and condition
on the event that the population will eventually hit the origin. Then
what will be the typical number of individuals that visit the origin,
and is there a limiting distribution for this number ?

\PA In the present work, we address a continuous version of the
previous problem, and so we deal with super-Brownian motion in
$\R^d$. We denote by $M_F(\R^d)$
the space of all finite measures in $\R^d$. We also
denote by $X=(X_t)_{t\geq 0}$ a
$d$-dimensional super-Brownian motion
with branching rate $\gamma$, that starts from
$\mu$ under the probability measure ${\PP}_\mu$, for every
$\mu\in M_F(\R^d)$. We refer to Perkins \cite{PE} for a
detailed presentation of super-Brownian motion. For every
$x\in\R^d$, we also denote by $\N_x$ the excursion measure
of super-Brownian motion
from $x$. We may and will assume that both $\PP_\mu$
and $\N_x$ are defined on the canonical space $C(\R_+,M_F(\R^d))$
of continuous functions from $\R_+$ into $M_F(\R^d)$ and that
$(X_t)_{t\geq 0}$ is the canonical process on this space. Recall from Theorem II.7.3 in Perkins \cite{PE}
that $X$ started at the Dirac measure $\delta_x$ can be constructed
from the atoms of a Poisson measure with intensity $\N_x$.

\PA The total
occupation measure of $X$ is the finite random measure on
$\R^d$ defined by
$${\mathcal Z}(A)=\int_0^\infty X_t(A)\,dt,$$
for every Borel subset $A$ of $\R^d$. We denote by $\mathcal R$
the topological support
of $\z$.
In dimension $d\geq 4$, points are polar, meaning that
$\N_x(0\in{\mathcal R})=0$ if $x\ne 0$, or equivalently
$\PP_\mu(0\in{\mathcal R})=0$ if $0$ does not belong to the
closed support of $\mu$. In dimension $d\leq 3$, we have
if $x\ne 0$,
\begin{equation}
\label{hitting0}
\N_x(0\in{\mathcal R})=\frac{8-2d}{\gamma}\,|x|^{-2}
\end{equation}
(see Theorem 1.3 in \cite{DIP} or Chapter VI in \cite{LG1}). It follows from the results in
Sugitani \cite{Su} that, again in dimension $d\leq 3$,
the measure $\z$ has a continuous density under
$\PP_{\delta_x}$ or under $\N_x$, for any $x\in\R^d$.
We write $(\ell^y,y\in\R^d)$ for this continuous
density.

\PA For every
$x\in\R^d$ and $r>0$, $B(x,r)$ denotes the open ball
centered at $x$ with radius $r$. To simplify notation,
we write $B_r=B(0,r)$ for the ball centered at
$0$ with radius $r$. By analogy with the
discrete problem mentioned above, we are
interested in the conditional distribution
of $\z(B_1)$ under $\PP_{\delta_x}(\cdot\mid \z(B_1)>0)$
when $|x|$ is large. As a simple consequence of
(\ref{hitting0}) and scaling, we have
when $d\leq 3$,
\begin{equation}
\label{probahitting}
\PP_{\delta_x}(\z(B_1)>0)
\sim
\N_x(\z(B_1)>0)
\sim
\frac{8-2d}{ \gamma}\,|x|^{-2}\qquad\hbox{as }|x|\to\infty.
\end{equation}
Here and later the notation $f(x)\sim g(x)$ as $|x|\to\infty$
means that the ratio $f(x)/g(x)$ tends to $1$ as $|x|\to\infty$.
On the other hand, when $d\geq 4$, it is proved in
\cite{DIP} that, as $|x|\to\infty$,
\begin{equation}
\label{probahitting2}
\PP_{\delta_x}(\z(B_1)>0)
\sim
\N_x(\z(B_1)>0)
\sim
\left\{
\begin{array}{ll}
{\displaystyle{2\over\gamma}} \,|x|^{-2}(\log|x|)^{-1}\qquad&\hbox{if }d=4,\\
\noalign{\medskip}
{\displaystyle{\kappa_d\over \gamma}}\,|x|^{2-d}&\hbox{if }d\geq 5,\\
\end{array}
\right.
\end{equation}
where $\kappa_d>0$ is a constant depending only on $d$.

\PA For $d\geq 3$, the Green function of $d$-dimensional Brownian
motion is 
$$G(x,y)=c_d|x-y|^{2-d},$$
where $c_d=(2\pi^{d/2})^{-1}\Gamma(\frac{d}{2}-1)$.
If $\mu\in M_F(\R^d)$ and $\varphi$ is a nonnegative measurable function 
on $\R^d$, we use the notation $\langle \mu,\varphi\rangle=\int \varphi\, d\mu$.
We can now state our main result.

\begin{theorem}
\label{main}
Let $\varphi$ be a bounded nonnegative measurable function supported on
$B_1$, and set
$\overline\varphi=\int \varphi(y)dy$.

\smallskip
\begin{enumerate}
\item[\rm(i)] If $d\leq 3$, the law of
$|x|^{d-4}\langle \z,\varphi\rangle$ under ${\PP}_{\delta_x}(\cdot\mid
{\mathcal
Z}(B_1)>0)$ converges as $|x|\to\infty$ towards the distribution
of $\overline\varphi\,\ell^0$ under $\N_{x_0}(\cdot \mid 0\in {\mathcal R})$,
where $x_0$ is an arbitrary point in $\R^d$ such that $|x_0|=1$.

\smallskip
\item[\rm(ii)] If $d=4$, the law of
$(\log|x|)^{-1}\langle \z,\varphi\rangle$ under ${\PP}_{\delta_x}(\cdot\mid
{\mathcal
Z}(B_1)>0)$ converges as $|x|\to\infty$ to an exponential distribution with mean 
$\gamma\overline\varphi/(4\pi^2)$.

\smallskip
\item[\rm(iii)] If $d\geq 5$, the law of
$\langle \z,\varphi\rangle$ under ${\PP}_{\delta_x}(\cdot\mid {\mathcal
Z}(B_1)>0)$ converges as $|x|\to\infty$ to the probability measure
$\mu_\varphi$ on $\R_+$ with moments $m_{p,\varphi}=\int
r^p\,\mu_\varphi(dr)$ given by
$$m_{1,\varphi}={c_d\over \kappa_d}\, \gamma\,\overline\varphi,$$
and for every $p\geq 2$,
$$m_{p,\varphi}={c_d\over \kappa_d}\, {\gamma^2\over 2}\,\sum_{j=1}^{p-1}
{p\choose j}
\int \N_z(\langle \z,\varphi\rangle^j)\N_z(\langle \z,\varphi\rangle^{p-j})dz.
$$
\end{enumerate}
\end{theorem}

\PA The scaling invariance properties of super-Brownian
motion allow us to restate Theorem \ref{main} in terms
of super-Brownian motion started with a fixed initial value and
the occupation measure of a small ball with radius $\varepsilon$
tending to $0$.
Part (i) of Theorem \ref{main} then becomes a straightforward
consequence of the fact that the measure $\mathcal Z$ has
a continuous density in dimension $d\leq 3$: See Lee \cite{Lee}
and Merle \cite{Merle} for more precise results along these
lines. On the other hand, the
proof of part (iii) is relatively easy from the method of moments
and known recursive formulas for the moments of the
random measure $\z$ under $\N_x$. For the sake of completeness,
we include proofs of the three cases in Theorem \ref{main},
but the most interesting part is really the
critical dimension $d=4$, where it is remarkable that an explicit limiting
distribution can be obtained.

  \PA Notice that
dimension $4$ is critical with respect to the polarity of points for
super-Brownian motion. Part (ii) of the theorem should therefore be
compared with classical limit theorems for additive functionals
of planar Brownian (note that $d=2$ is the critical dimension
for polarity of points for ordinary Brownian motion). The celebrated
Kallianpur-Robbins law states that the time spent by planar Brownian
motion in a
bounded set before time $t$ behaves as $t\to\infty$
like $\log t$ times an exponential variable (see e.g. section 7.17
in It\^o and McKean \cite{IM}). The
Kallianpur-Robbins law can be derived by ``conceptual proofs'' which
explain the occurence of the exponential
distribution. Our initial approach to part (ii) was based on
a similar conceptual argument based on the Brownian snake
approach to super-Brownian motion. Since it seems delicate to
make this argument completely rigorous, we rely below on
a careful analysis of the moments of $\langle \z,\varphi\rangle$.

\PA Let us finally comment on the branching random walk problem
discussed at the beginning of this introduction. Although we 
do not consider this problem here, it is very likely that a result
analogous to Theorem \ref{main} holds in this discrete setting, just
replacing $\langle\z,\varphi\rangle$ with the number of particles 
that hit the origin. In particular, the limiting distributions 
obtained in (i) and (ii) of Theorem \ref{main} should also appear
in the discrete setting.

\section{Preliminary remarks}

Let us briefly recall some basic facts about
super-Brownian motion and its excursion measures. If $x\in\R^d$
and if
$${\mathcal N}=\sum_{i\in I} \delta_{\omega_i}$$
is a Poisson point measure on $C(\R_+,M_F(\R^d))$ with intensity 
$\N_x(\cdot)$, then the mea\-sure-valued process $Y$ defined
by 
\begin{eqnarray*}
&&Y_0=\delta_x,\\
&&Y_t=\sum_{i\in I} X_t(\omega_i)\ ,\hbox{ for every }t>0,
\end{eqnarray*}
has law $\PP_{\delta_x}$ (see Theorem II.7.3 in \cite{PE}). 

\PA We can use this Poisson decomposition to observe  that it is enough to prove
Theorem \ref{main} with the conditional measure 
${\PP}_{\delta_x}(\cdot\mid {\mathcal
Z}(B_1)>0)$ replaced by $\N_x(\cdot\mid {\mathcal
Z}(B_1)>0)$. Indeed,
write $M=\#\{i\in I: {\mathcal R}(\omega_i)\cap B_1\ne \emptyset\}$. Then, 
$M$ is Poisson with
parameter $\N_x(\z(B_1)>0)$, and $\{M\geq 1\}$ is the event that the
range of $Y$ hits $B_1$. Furthermore,
the preceding Poisson decomposition just shows that
the law of $\langle \z,\varphi\rangle$ under ${\PP}_{\delta_x}$ coincides
with the
law of $Z^1+\cdots+Z^{M}$, where conditionally given $M$, the variables $Z^1,Z^2,\ldots$
are independent and distributed according to the law
of $\langle \z,\varphi\rangle$ under $\N_x(\cdot \mid \z(B_1)>0)$. Since
$P(M=1\mid
M\geq 1)$
tends to $1$ as $|x|\to\infty$ (by the estimates (\ref{probahitting})
and (\ref{probahitting2})), we see that the law of $\langle
\z,\varphi\rangle$ (or the law of $f(x)\langle
\z,\varphi\rangle$ for any deterministic function $f$) under
${\PP}_{\delta_x}(\cdot\mid {\mathcal
Z}(B_1)>0)$ will be arbitrarily close to the law of the same variable
under $\N_x(\cdot\mid {\mathcal
Z}(B_1)>0)$ when $|x|$ is large, which is what we wanted. Note that this
argument is valid in any dimension.

\PA Let us also discuss the dependence of our results on the branching rate $\gamma$.
If $(Y_t)_{t\geq 0}$ is a super-Brownian motion with branching rate $\gamma$
started at $\mu$, and $\lambda>0$, then $(\lambda Y_t)_{t\geq 0}$ is a super-Brownian motion with branching rate 
$\lambda\gamma$
started at $\lambda\mu$. A similar property then holds for excursion measures.
Write $\N_x^{(\gamma)}$ instead of $\N_x$ to emphasize the dependence on $\gamma$.
Then the ``law'' of $(\lambda X_t)_{t\geq 0}$ under $\N^{(\gamma)}_x$
is $\lambda^{-1}\N_x^{(\lambda\gamma)}$. Thanks to these observations, it will
be enough to prove Theorem \ref{main} for one particular value of $\gamma$.

\PA In what follows, we take $\gamma=2$, as this will simplify certain
formulas. For any nonnegative measurable function $\varphi$
on $\R^d$, the moments
of $\langle\z,\varphi\rangle$ are determined by induction
by the formulas
\begin{equation}
\label{momentZ1}
\N_x(\langle\z,\varphi\rangle)=\int_{\R^d} G(x,y)\varphi(y)\,dy
\end{equation}
and, for every $p\geq 2$,
\begin{equation}
\label{momentZp}
\N_x(\langle\z,\varphi\rangle^p)
=\sum_{j=1}^{p-1}
{p\choose j}
\int_{\R^d} G(x,z)\,\N_z(\langle \z,\varphi\rangle^j)\N_z(\langle
\z,\varphi\rangle^{p-j})dz.
\end{equation}
See e.g. formula (16.2.3) in \cite{LG2}, and note that the extra 
factor $2$ there is due to the fact that the 
the Brownian snake approach gives $\gamma=4$.

\section{Low dimensions}

\PA In this section, we prove part (i) of Theorem \ref{main}.
Let $\varepsilon>0$, and set
$\varphi_\varepsilon(x)=\varphi(x/\varepsilon)$. By scaling,
the law of
$\langle\z,\varphi\rangle$ under $\N_x(\cdot\mid {\mathcal
Z}(B_1)>0)$ coincides with the law of
$\varepsilon^{-4}\langle\z,\varphi_\varepsilon\rangle$ under
$\N_{\varepsilon x}(\cdot\mid {\mathcal
Z}(B_\varepsilon)>0)$. Taking $\varepsilon=|x|^{-1}$, we see that the
proof of part (i)
reduces to checking that the law of
$|x|^{d}\langle\z,\varphi_{|x|^{-1}}\rangle$
under $\N_{x/|x|}(\cdot\mid {\mathcal
Z}(B_{|x|^{-1}})>0)$ converges to the distribution of
$\overline\varphi\ell^0$
under $\N_{x_0}(\cdot\mid0\in{\mathcal R})$. 

\PA However, as $|x|\to\infty$,
$$\N_{x/|x|}({\mathcal
Z}(B_{|x|^{-1}})>0)=\N_{x_0}(\z(B_{|x|^{-1}})>0)\longrightarrow
\N_{x_0}(0\in{\mathcal R})=4-d.$$
On the other hand, since
$$|x|^{d}\langle\z,\varphi_{|x|^{-1}}\rangle
=|x|^d \int dy\,\ell^y\,\varphi_{|x|^{-1}}(y)=\int dy\,\ell^{y/|x|}\,\varphi(y)$$
the continuity of the local times $\ell^y$ implies that, for every
$\delta>0$,
$$\N_{x/|x|}\Big(
\Big||x|^{d}\langle\z,\varphi_{|x|^{-1}}\rangle-\overline\varphi\,\ell^0\Big|>\delta\Big)
\leq \N_{x/|x|}\Big(\sup_{y\in
B_{|x|^{-1}}}|\ell^y-\ell^0|>\frac{\delta}{\overline\varphi}\Big)
\longrightarrow 0$$
as $|x|\to\infty$. By rotational invariance, the law of $\ell^0$ under $\N_{x/|x|}$
coincides with the law of the same variable under $\N_{x_0}$. Part (i)
of Theorem \ref{main} now follows from the preceding observations.

\section{High dimensions}

\PA We now turn to part (iii) of Theorem \ref{main}
and so we suppose that $d\geq 5$. As noticed earlier, we may replace
${\PP}_{\delta_x}(\cdot\mid {\mathcal
Z}(B_1)>0)$ by $\N_x(\cdot\mid {\mathcal
Z}(B_1)>0)$. 

\PA Without loss of generality, we assume in this part that $\varphi\leq
1$.

\begin{lemma}
\label{tech}
There exists a finite constant $K_d$ depending only
on $d$, such that, for every $x\in\R^d$ and $p\geq 1$,
$$\N_x(\langle\z,\varphi\rangle^p)\leq K_d^p\,p!\,(|x|^{2-d}\wedge 1).$$
\end{lemma}

\noindent{\bf Proof.} Obviously, it is enough to consider
the case when $\varphi={\bf 1}_{B_1}$. From (\ref{momentZ1}), one immediately
verifies that
$$\N_x(\z(B_1))\leq C_{1,d}\,(|x|^{2-d}\wedge 1)$$
for some constant $C_{1,d}$ depending only on $d$. Straightforward
estimates give the existence of a constant $a_d$ such that,
for every $x\in\R^d$,
$$\int G(x,z)\,(|z|^{2-d}\wedge 1)^2\,dz \leq a_d(|x|^{2-d}\wedge 1).$$
We then claim that for every integer $p\geq 1$,
\begin{equation}
\label{tech2}
\N_x(\z(B_1)^p)\leq C_{p,d}\,p!\,(|x|^{2-d}\wedge 1)
\end{equation}
where the constants $C_{p,d}$, $p\geq 2$ are determined by induction by
\begin{equation}
\label{tech3}
C_{p,d}=a_d\sum_{j=1}^{p-1} C_{j,d}C_{p-j,d}.
\end{equation}
Indeed, let $k\geq 2$ and suppose that (\ref{tech2}) holds for every
$p\in\{1,\ldots,k-1\}$. From (\ref{momentZp}), we get
$$\N_x(\z(B_1)^k)
\leq k!\sum_{j=1}^{k-1} C_{j,d}C_{k-j,d}\int G(x,z)\,(|z|^{2-d}\wedge
1)^2\,dz,$$
and our choice of $a_d$ shows that (\ref{tech2}) also holds for $p=k$.
We have thus proved our claim (\ref{tech2}) for every $p\geq 1$.

\PA From (\ref{tech3}) it is an elementary exercise to verify that
$C_{p,d}\leq K_d^p$ for some constant $K_d$ depending only on $d$. This
completes the proof. \hfill$\square$

\smallskip
\PA Let us now prove that for every $p\geq 1$, $\N_x(\langle
\z,\varphi\rangle^p
\mid \z(B_1)>0)$ converges as $|x|\to\infty$ to $m_{p,\varphi}$. If
$p=1$, this is an immediate consequence of (\ref{probahitting2})
and (\ref{momentZ1}). If $p\geq 2$, we write
$${\N_x(\langle\z,\varphi\rangle^p)\over |x|^{2-d}}
=\sum_{j=1}^{p-1}
{p\choose j}
\int {G(x,z)\over |x|^{2-d}}\,\N_z(\langle \z,\varphi\rangle^j)\N_z(\langle
\z,\varphi\rangle^{p-j})dz
$$
and we can use the bounds of Lemma \ref{tech}, and the property $\int (|z|^{2-d}\wedge 1)^2dz<\infty$, in order
to get
$$\lim_{|x|\to\infty}
{\N_x(\langle\z,\varphi\rangle^p)\over |x|^{2-d}}
=c_d
\sum_{j=1}^{p-1}
{p\choose j}
\int \N_z(\langle \z,\varphi\rangle^j)\N_z(\langle
\z,\varphi\rangle^{p-j})dz.$$
The convergence of $\N_x(\langle \z,\varphi\rangle^p
\mid \z(B_1)>0)$ towards $m_{p,\varphi}$ now follows from
(\ref{probahitting2}).

\PA Finally, Lemma \ref{tech} and (\ref{probahitting2}) also imply 
that  any limit distribution of the
laws of $\langle \z,\varphi\rangle$ under
$\N_x(\cdot
\mid \z(B_1)>0)$ is characterized by its moments. Part (iii)
of Theorem \ref{main} now follows as a standard application of
the method of moments.

\section{The critical dimension}

In this section, we consider the critical dimension $d=4$. Recall that 
in that case $G(x,y)=(2\pi^2)^{-1}|y-x|^{-2}$. As in the previous sections, we take $\gamma=2$. We start by stating
two lemmas.

\begin{lemma}\label{estipositive}
Let $x \in \R^4\setminus\{0\}$. Then,
\begin{equation*}
\mathbb{N}_{x}[\z(B_\varepsilon)>0 ]
\underset{\varepsilon \to 0}{\sim}
|x|^{-2}\left(\log \frac{1}{\varepsilon} \right)^{-1}.
\end{equation*}
\end{lemma}

\begin{lemma}\label{estimomentp}
Let $x \in \R^4\setminus\{0\}$, and $p \geq 1$. Let
$\varphi$ be a bounded nonnegative measurable function on $B_1$, and for every
$\varepsilon>0$, put $\varphi_\varepsilon(y)=\varphi(y/\varepsilon)$. Then,
\begin{equation*}
\mathbb{N}_{x}[ \langle \z,\varphi_\varepsilon\rangle)^{p}
] \underset{\varepsilon \to 0}{\sim}
p!\,\left(\frac{\overline\varphi}{2\pi^2}\right)^p \,|x|^{-2}\,\varepsilon^{4p} \left(
\log \frac{1}{\varepsilon} \right)^{p-1},
\end{equation*}
uniformly when $x$ varies over a compact subset of $\R^4\setminus\{0\}$.
\end{lemma}

\PA Let us explain how part (ii) of Theorem \ref{main}
follows from these two lemmas. Notice that the estimate of Lemma 
\ref{estipositive} also holds  
uniformly when $x$ varies over a compact subset of $\R^4\setminus\{0\}$,
by scaling and rotational invariance. Combining the results of the lemmas gives
$$\mathbb{N}_{x}\left[ \left(
\frac{\langle
\z,\varphi_\varepsilon\rangle}{\varepsilon^{4}\log(\frac{1}{\varepsilon})}
\right)^{p} \bigg|\z(B_\varepsilon)>0 \right]
\underset{\varepsilon \to 0}{\sim} p!\left(\frac{\overline\varphi}{2\pi^2}\right)^p,$$
uniformly when $x$ varies over a compact subset of $\R^4\setminus\{0\}$. By scaling, for any $x\in\R^4$ with
$|x|>1$ the law of
$\langle \z,\varphi\rangle$ under $\N_x(\cdot\mid \z(B_1)>0)$
coincides with the law of
$|x|^4\langle \z,\varphi_{1/|x|}\rangle$ under $\N_{x/|x|}(\cdot\mid
\z(B_{1/|x|})>0)$.
Hence, we deduce from the preceding display that we have also
$$\mathbb{N}_{x}\left[ \left(
\frac{\langle \z,\varphi\rangle}{\log |x|}
\right)^{p} \bigg|\z(B_1)>0 \right]
\underset{|x|\to\infty}{\sim} p!\left(\frac{\overline\varphi}{2\pi^2}\right)^p.$$
The statement in part (ii) of Theorem \ref{main} now follows from
an application of the method of moments. 
\PA It remains to prove
Lemma \ref{estipositive} and Lemma \ref{estimomentp}.

\subsection{Proof of Lemma \ref{estipositive}}
From well-known connections between super-Brownian motion and 
partial differential equations (see e.g. Chapter VI in \cite{LG1}), the function
$u_\varepsilon(x)=\N_x[\z(B_\varepsilon)>0]$
defined for $|x|>\varepsilon$
solves the singular boundary problem
\begin{eqnarray*}
&&\Delta u =2\,u^2,\qquad \hbox{in the domain }\{|x|>\varepsilon\}\\
&&u(x)\longrightarrow \infty, \quad\;\hbox{as }|x|\to\varepsilon^+\\
&&u(x)\longrightarrow 0, \qquad\hbox{as }|x|\to \infty.
\end{eqnarray*}
As a consequence of a lemma due to Iscoe (see Lemma 3.4 in \cite{DIP}), for every
$x\in\R^4\backslash \{0\}$, we have $u_\varepsilon(x)\sim
|x|^{-2}(\log(1/\varepsilon))^{-1}$
as $\varepsilon\to 0$. Lemma \ref{estipositive} follows.

\newpage

\subsection{Proof of Lemma \ref{estimomentp}}

\subsubsection{Lower bound}
Let us introduce the space of functions
$$\mathcal{F} :\ = \Big\{ f: \mathbb{R}_{+} \rightarrow
\mathbb{R} :
   \lim_{\varepsilon \to 0} f(\varepsilon) =0  \mbox{ and }
   \lim_{\varepsilon \to 0} \frac{f(\varepsilon)}{\varepsilon}
   =\infty \Big \}. $$

{\bf Claim}. {\it For every integer $p\geq 1$, for every $f\in\mathcal{F}$,
for every $\beta >0$, there exists $\varepsilon_0>0$ such that,
for every $\varepsilon\in(0,\varepsilon_0)$,}
\begin{equation}
\label{lowerestimate}
\inf_{y \notin B_{f(\varepsilon)}}
 |y|^{2} \left(\log\left(\frac{|y|}{\varepsilon}
\right)\right)^{1-p}
\mathbb{N}_{y}\left[  \langle\z,\varphi_\varepsilon\rangle^{p} \right]
   \geq (1-\beta)\,p! \left(\frac{\overline\varphi}{2\pi^2}\right)^p\,\varepsilon^{4p}.
\end{equation}

\PA We prove the claim by induction on $p$.   Let us first consider the
case $p=1$.
We fix $f \in \mathcal{F}$. Using (\ref{momentZ1}),
for $\varepsilon>0$ and $y\in\R^4$ such that $|y| > f(\varepsilon)$, we have
\begin{equation} \label{proofH1}
\mathbb{N}_{y}\left[ \langle \z,\varphi_\varepsilon\rangle \right]=
\int_{B_\varepsilon} dz
\varphi_\varepsilon(z)
G(y,z) \geq
\varepsilon^{4}\overline\varphi  \inf_{z \in
B_\varepsilon}G(y,z) .
\end{equation}
Since $\lim_{\varepsilon \to 0 } (f(\varepsilon)/\varepsilon)=
\infty $, we
see that
   $$\inf_{y \notin B_{f(\varepsilon)}} \Big(|y|^{2}\inf\big\{
G(y,z),z \in B_\varepsilon\big\}\Big)
\underset{\varepsilon \to 0}{\rightarrow} \frac{1}{2\pi^{2}}. $$
   We thus deduce from (\ref{proofH1}) that for $\varepsilon$ small enough,
$$\inf_{y \notin B_{f(\varepsilon)}}
 |y|^{2}
\mathbb{N}_{y}\left[ \langle \z,\varphi_\varepsilon\rangle \right] \geq
\frac{1-\beta}{2\pi^2}\,\overline\varphi\,\varepsilon^4,$$
which gives our claim for $p=1$.
\vskip0.1cm \PA
Let $p \geq 2$ and suppose that the claim holds up to order $p-1$.
 Fix $f \in \mathcal{F}$ and $\beta\in(0,1)$. Let $\beta'\in(0,1)$
be such that $(1-\beta')^4=1-\beta$, and  let $C>0$ be
such that $(1+C^{-1})^{-2} = 1- \beta'.$ Introduce
   the function $\hat{f}$ defined by
\begin{equation}
\label{defifhat}
\hat{f}(\varepsilon) = \varepsilon
\log\left(\frac{f(\varepsilon)}{\varepsilon}\right).
\end{equation}
Clearly, $\hat{f} \in \mathcal{F}$. Furthermore, we have
\begin{equation} \label{limitffhat}
   \lim_{\varepsilon \to 0}
\frac{\log(\hat{f}(\varepsilon)/\varepsilon)}{\log(
f(\varepsilon)/\varepsilon)} =0.
\end{equation}
Using (\ref{momentZp}),
we obtain, for any $y \notin B_{f(\varepsilon)}$
$$\mathbb{N}_{y}\left[  \langle \z,\varphi_\varepsilon\rangle ^{p}
\right]
    \geq \sum_{j=1}^{p-1}
{p \choose j}
\int_{B_{|y|/C} \setminus B_{\hat{f}(\varepsilon)}} dz\,
G(y,z)\, \mathbb{N}_{z}\left[\langle \z,\varphi_\varepsilon\rangle^{j} \right]
\mathbb{N}_{z}\left[ \langle \z,\varphi_\varepsilon\rangle^{p-j}\right].
$$
Using the induction hypothesis,
we get, provided $\varepsilon$ is small enough,
\begin{eqnarray*}
&&\mathbb{N}_{y}\left[  \langle \z,\varphi_\varepsilon\rangle
^{p} \right]\\
&&\quad\geq(1-\beta')^2p!
\left(\frac{\overline\varphi}{2\pi^2}\right)^p 
(p-1)\varepsilon^{4p}
\int_{B_{|y|/C}\setminus
B_{\hat{f}(\varepsilon)}}
   dz
G(y,z)   \frac{1}{|z|^{4}} \left(
\log\frac{|z|}{\varepsilon} \right)^{p-2}.
\end{eqnarray*}
From the definition of $C$, for any $z \in B_{|y|/C}$, we have
$G(y,z)\geq
(1-\beta') G(0,y)$.
It follows that
\begin{eqnarray*}
&&(p-1)\int_{B_{|y|/C}\setminus
B_{\hat{f}(\varepsilon)}} dz
G(y,z)   \frac{1}{|z|^{4}} \left(
\log\frac{|z|}{\varepsilon} \right)^{p-2}\\
&&\qquad\geq2\pi^2 (1-\beta')\,G(0,y)(p-1)
\int_{\hat{f}(\varepsilon)}^{|y|/C}
     \frac{dr}{r} \left(
\log \frac{r}{\varepsilon} \right)^{p-2}\\
&&\qquad=(1-\beta')|y|^{-2}\left(\left(\log 
\frac{|y|}{C\varepsilon}\right)^{p-1} -
\left(\log \frac{\hat{f}(\varepsilon)}{\varepsilon} \right)^{p-1}\right).
\end{eqnarray*}
Moreover, using the property $f\in{\mathcal F}$ and (\ref{limitffhat}), we see that, if $\varepsilon$ is
sufficiently small, for any $y \notin B_{f(\varepsilon)}$
$$ \left(\log  \frac{|y|}{C\varepsilon}\right)^{p-1} -
\left(\log \frac{\hat{f}(\varepsilon)}{\varepsilon} \right)^{p-1}
\geq (1-\beta') \left( \log\frac{|y|}{\varepsilon}  \right)^{p-1}.$$
From the preceding bounds, we get that, if $\varepsilon$ is sufficiently
small,
$$\inf_{y \notin B_{f(\varepsilon)}}
|y|^{2} \left( \log\frac{|y|}{\varepsilon} \right)^{1-p}
\mathbb{N}_{y}\left[ \langle \z,\varphi_\varepsilon\rangle
^{p} \right] \geq   (1-\beta')^{4} p!
\left(\frac{\overline\varphi}{2\pi^2}\right)^p\,\varepsilon^{4p},
$$
which is our claim at order $p$.

\subsubsection{Upper bound}
Without loss of generality we assume that $\varphi\leq 1$.
We need to get
upper bounds on $\mathbb{N}_{y}\left[  
\langle \z,\varphi_\varepsilon\rangle^{p} \right]$ for $y$ belonging to different
subsets of $\mathbb{R}^{4}$.

\PA   We will prove that, for every $p\geq 1$, for every $f\in{\mathcal F}$  and every $\beta\in(0,1)$
the following bounds hold for $\varepsilon>0$ sufficiently small:
\begin{eqnarray*}
\left(\mathfrak{H}_{p}^{1}\right)  \!&&\!
\sup_{|y| \leq
4\varepsilon}
\mathbb{N}_{y}\left[
\langle \z,\varphi_\varepsilon\rangle^{p} \right] \leq p!\,
\varepsilon^{4p-2} \left( \log\frac{f(\varepsilon)}{\varepsilon}
\right)^{p-1}  , \\
\left(\mathfrak{H}_{p}^{2}\right) \!&&\!
\sup_{4\varepsilon \leq |y| \leq
f(\varepsilon)}\!\!\!
|y|^{2} \mathbb{N}_{y}\left[ \langle \z,\varphi_\varepsilon\rangle^{p}
\right] \leq p!\,\varepsilon^{4p}\left(\log 
\frac{f(\varepsilon)}{\varepsilon} 
\right)^{p-1}  , \\
   \left(\mathfrak{H}_{p}^{3}\right) \!&&\!  \sup_{|y| \geq f(\varepsilon)}
 |y|^{2} \left( \log
\frac{|y|}{\varepsilon}\right)^{1-p}
\mathbb{N}_{y}\left[ \langle \z,\varphi_\varepsilon\rangle^{p}
\right] \leq \left(\frac{\overline\varphi+\beta}{2\pi^2}\right)^p p!\;\varepsilon^{4p} .
   \end{eqnarray*}
Only $(\mathfrak{H}_{p}^{3})$ is needed
in our proof of Lemma  \ref{estimomentp}. However, we will proceed by induction 
on $p$ to get $(\mathfrak{H}_{p}^{3})$, and we will use $(\mathfrak{H}_{p}^{1})$
and $(\mathfrak{H}_{p}^{2})$ in our induction argument. The bounds $(\mathfrak{H}_{p}^{1})$
and $(\mathfrak{H}_{p}^{2})$ are not sharp, but they will be sufficient for our purposes.
 Notice that 
$(\overline \varphi +\beta)/(2\pi^2)<1/3$ because $\varphi\leq 1$ and $\beta<1$.
\vskip0.1cm

\PA   We first note that when $p=1$ the bounds $(\mathfrak{H}_{1}^1)$, $(\mathfrak{H}_{1}^2)$
and $(\mathfrak{H}_{1}^3)$ are easy consequences of
(\ref{momentZ1}). Let $p \geq 2$ and
assume that $(\mathfrak{H}_{k}^1)$, $(\mathfrak{H}_{k}^2)$
and $(\mathfrak{H}_{k}^3)$ hold for every $1\leq k \leq
p-1$, for any choice of $\beta$ and $f$. Let us fix $f \in \mathcal{F}$ and $\beta_0\in(0,1)$.
In our induction argument we will use $(\mathfrak{H}_{k}^3)$, for $1\leq k\leq p-1$, with
$\beta\in(0,1)$ chosen small enough so that
$$(1+\beta)^2(\overline \varphi+\beta)^p <
(\overline\varphi+\beta_0)^p.$$

\PA For every $j\in\{1,\ldots,p-1\}$, every $y\in\R^4$
and every Borel subset $A$ of $\R^4$, we set
$$I^\varepsilon_{p,j}(A,y):=\frac{1}{j! (p-j)!}\int_{A} dz G(y,z)
\mathbb{N}_{z}\left[ \langle \z,\varphi_\varepsilon\rangle^{j} \right]
\mathbb{N}_{z}\left[ \langle \z,\varphi_\varepsilon\rangle^{p-j}\right].$$
We also set $I^\varepsilon_{p,j}(y)=I^\varepsilon_{p,j}(\R^4,y)$.

\PA We first verify $(\mathfrak{H}_{p}^{1})$, and so we assume that
$|y|\leq 4\varepsilon$. We fix $j\in\{1,\ldots,p-1\}$ and we split the
the integral in $I^\varepsilon_{p,j}(y)$ into
three parts corresponding to the sets
   $$A_{1}^{(1)}=B_{ 8\varepsilon}, \ A_{2}^{(1)} = B_{ f(\varepsilon)}
\setminus
B_{8\varepsilon}, \
   A_{3}^{(1)}= \mathbb{R}^{4} \setminus B_{
   f(\varepsilon)}.$$
From (\ref{momentZp}), we have
\begin{equation} \label{decoupypetit}
\mathbb{N}_{y} \left[ \langle \z,\varphi_\varepsilon\rangle^{p} \right] = p! \sum_{j=1}^{p-1} 
\left( I^\varepsilon_{p,j}(A_{1}^{(1)},y) + I^\varepsilon_{p,j}(A_{2}^{(1)},y) +
I^\varepsilon_{p,j}(A_{3}^{(1)},y) \right).
\end{equation}
If $\varepsilon$ is small enough, we deduce from the bounds
$(\mathfrak{H}_{k}^{1})$ and $(\mathfrak{H}_{k}^{2})$, with ${1 \leq k \leq p-1}$, that
$$I^\varepsilon_{p,j}(A_{1}^{(1)},y) \leq 
\varepsilon^{4p-4} \left( \log\frac{f(\varepsilon)}{\varepsilon}
\right)^{p-2} \int_{B_{8\varepsilon}}
\frac{dz}{2\pi^{2}|z-y|^{2}}.$$
It follows that
$$ \lim_{\varepsilon \to 0} \left(\sup_{|y|\leq 4 \varepsilon}
\varepsilon^{2-4p} \left( \log \frac{f(\varepsilon)}{\varepsilon}
\right)^{1-p} I^\varepsilon_{p,j}(A_{1}^{(1)},y) \right)=0.$$
\vskip0.1cm
If $z \in A_{2}^{(1)}\cup A_3^{(1)}$, we have
$G(y,z) \leq 4 G(0,z)$.
Using $(\mathfrak{H}_{k}^{2})$ with $1 \leq k \leq p-1$, we obtain,
if $\varepsilon$ is small enough,
\begin{eqnarray*}
I^\varepsilon_{p,j}(A_{2}^{(1)},y) & \leq & 4\,
\varepsilon^{4p} \left( \log \frac{f(\varepsilon)}{\varepsilon}
\right)^{p-2} \int_{8\varepsilon}^{f(\varepsilon)} \frac{dr}{r^{3}},
\end{eqnarray*}
and thus
$$ \lim_{\varepsilon \to 0}\left(
\sup_{|y| \leq 4\varepsilon}\varepsilon^{2-4p}
\left( \log \frac{f(\varepsilon)}{\varepsilon} \right)^{1-p}
I^\varepsilon_{p,j}(A_{2}^{(1)},y) \right)=0. $$
\vskip0.1cm
   Using finally $(\mathfrak{H}_{k}^{3})$ with $1\leq k \leq p-1$,
we obtain, for $\varepsilon$ sufficiently small,
   $$I^\varepsilon_{p,j}(A_{3}^{(1)},y) \leq 4\,
\varepsilon^{4p} \int_{f(\varepsilon)}^{\infty}
\frac{dr}{r^{3}} \left( \log \frac{r}{\varepsilon}
\right)^{p-2}.$$ 
Since
$$\int_{f(\varepsilon)}^{\infty}
\frac{dr}{r^{3}} \left( \log \frac{r}{\varepsilon}
\right)^{p-2}\underset{\varepsilon \to 0}{\sim} 
\frac{1}{2f(\varepsilon)^2}\,\left( \log \frac{f(\varepsilon)}{\varepsilon}
\right)^{p-2}$$
we get
$$ \lim_{\varepsilon \to 0} \left(\sup_{|y| \leq 4 \varepsilon}
\varepsilon^{2-4p} \left( \log
\frac{f(\varepsilon)}{\varepsilon}
\right)^{1-p} I^\varepsilon_{p,j}(A_{3}^{(1)},y)\right)=0.$$
   \vskip0.1cm
Combining the estimates we
obtained for $I^\varepsilon_{p,j}(A_{1}^{(1)},y),I^\varepsilon_{p,j}(A_{2}^{(1)},y)$ and
$I^\varepsilon_{p,j}(A_{3}^{(1)},y)$,
we arrive at
$$ \lim_{\varepsilon \to 0}\left( \sup_{|y| \leq 4 \varepsilon}
\varepsilon^{2-4p} \left( \log
\frac{f(\varepsilon)}{\varepsilon}
\right)^{1-p} I^\varepsilon_{p,j}(y)\right) =0.$$
From (\ref{decoupypetit}), we obtain that $(\mathfrak{H}_{p}^{1})$ holds.

\vskip0.2cm
\PA We now turn to the proof of $(\mathfrak{H}_{p}^{2})$, and so we assume
that $4 \varepsilon \leq |y| \leq f(\varepsilon)$. Again we fix $j\in\{1,\ldots,p-1\}$.
We split the integral in
$I^\varepsilon_{p,j}(y)$ into five parts corresponding to the sets
\begin{itemize}
\item
$A_{1}^{(2)}=B_{ 2\varepsilon}$,
\item
   $A_{2}^{(2)} = B(y,|y|/2)$,
\item
$A_{3}^{(2)}= B_{\hat{f}(\varepsilon)}\setminus \left(
B_{ 2\varepsilon} \cup  B(y,|y|/2)     \right)$,
\item
   $A_{4}^{(2)}= B_{2f(\varepsilon)} \setminus \left( B_{
   \hat{f}(\varepsilon)} \cup B(y, |y|/2) \right)$,
\item
$A_{5}^{(2)} = \mathbb{R}^{4} \setminus B_{2f(\varepsilon)}$,
\end{itemize}
where $\hat{f}(\varepsilon)=\varepsilon  \log
(f(\varepsilon)/\varepsilon)$ as in (\ref{defifhat}). We have thus
\begin{equation} \label{decoupyinter}
\mathbb{N}_{y} \left[ \langle \z,\varphi_\varepsilon\rangle^{p} \right] = p!\sum_{j=1}^{p-1} 
\sum_{i=1}^{5} I^\varepsilon_{p,j}(A_{i}^{(2)},y).
\end{equation}
Notice first that if $z \in A_{1}^{(2)}$, we have  $|z| \leq |y|/2$ so that
$G(z-y) \leq 4G(y)$.
Using $(\mathfrak{H}_{k}^{1})$ with $1 \leq k \leq p-1$, we obtain, provided
$\varepsilon$ is small enough
   $$ I^\varepsilon_{p,j}(A_{1}^{(2)},y) \leq 
\varepsilon^{4p-4} \left( \log 
\frac{f(\varepsilon)}{\varepsilon} 
   \right)^{p-2} \frac{2}{\pi^{2}|y|^{2}} \int_{\{|z| \leq 2 \varepsilon\}}
dz, $$
so that
$$ \lim_{\varepsilon \to 0} \left(\sup_{4\varepsilon < |y| < f(\varepsilon)}
\!\! \varepsilon^{-4p} |y|^{2} \left( \log 
\frac{f(\varepsilon)}{\varepsilon} 
   \right)^{1-p} I^\varepsilon_{p,j}(A_{1}^{(2)},y)\right) =0. $$
If $z \in A_{2}^{(2)}$, using the bound
$|z|^{-2} \leq 4|y|^{-2} $, we deduce from $(\mathfrak{H}_{k}^{1})$ and
$(\mathfrak{H}_{k}^{2})$ for $1 \leq k \leq p-1$ that for sufficiently small
$\varepsilon$,
$$  I^\varepsilon_{p,j}(A_{2}^{(2)},y) \leq 
\varepsilon^{4p} \left( \log
\frac{f(\varepsilon)}{\varepsilon} 
   \right)^{p-2} \frac{16^3}{|y|^{4}}
\int_{B(y,|y|/2)} G(y,z) dz. $$
It follows that
$$ \lim_{\varepsilon \to 0} \left(\sup_{4\varepsilon < |y| < f(\varepsilon)}
\!\! \varepsilon^{-4p} |y|^{2}\left( \log 
\frac{f(\varepsilon)}{\varepsilon} 
   \right)^{1-p} I^\varepsilon_{p,j}(A_{2}^{(2)},y) \right)=0. $$
If $z \in A_{3}^{(2)}$, $G(y,z)
\leq 4G(0,y) $.
Since $\hat{f} \in \mathcal{F}$, we can use
$(\mathfrak{H}_{k}^{1})$ and   $(\mathfrak{H}_{k}^{2})$ with $1 \leq k \leq p-1$
to get that, for $ \varepsilon$ small enough,
$$ I^\varepsilon_{p,j}(A_{3}^{(2)},y)
\leq \varepsilon^{4p} \left( \log
\frac{\hat{f}(\varepsilon)}{\varepsilon}
\right)^{p-2}
\frac{4\times 16^2}{|y|^{2}} \int_{\varepsilon}^{\hat{f}(\varepsilon)} r^{-1}
dr. $$
It then follows from (\ref{limitffhat}) that
$$ \lim_{\varepsilon \to 0} \left(\sup_{4\varepsilon < |y| < f(\varepsilon)}
\varepsilon^{-4p} |y|^{2} \left( \log 
\frac{f(\varepsilon)}{\varepsilon} 
   \right)^{1-p} I^\varepsilon_{p,j}(A_{3}^{(2)},y)\right) =0 . $$
If $z \in A_{4}^{(3)}$, we still have $G(y,z)
\leq 4G(0,y) $. Again, $\hat{f} \in \mathcal{F}$,
and we can use
$(\mathfrak{H}_{k}^{3})$ with $1 \leq k \leq p-1$, recalling that
$(\overline \varphi +\beta)/(2\pi^2)<1/3$, to obtain
for $\varepsilon$ small
\begin{eqnarray*} I^\varepsilon_{p,j}(A_{4}^{(2)},y) &\leq &3^{-p}\,
\varepsilon^{4p}
   \frac{4}{|y|^{2}}
\int_{\hat{f}(\varepsilon)}^{2f(\varepsilon)}
\frac{dr}{r} \left(  \log \frac{r}{\varepsilon} \right)^{p-2}
\\ & = &\frac{4\times 3^{-p}}{p-1}\;
\varepsilon^{4p} \frac{1}{|y|^{2}}
   \left(  \log \frac{2f(\varepsilon)}{\hat f(\varepsilon)} 
\right)^{p-1}.
\end{eqnarray*}
It follows that
$$ \limsup_{\varepsilon \to 0} \left(\sup_{4\varepsilon < |y| < f(\varepsilon)}
\!\! \varepsilon^{-4p} |y|^{2}\left( \log 
\frac{f(\varepsilon)}{\varepsilon} 
   \right)^{1-p} I^\varepsilon_{p,j}(A_{4}^{(2)},y)\right) 
\leq \frac{4\times 3^{-p}}{p-1}<\frac{1}{p-1}. $$
Finally, if $|z|\geq 2f(\varepsilon)$, we have $G(y,z)\leq 4(2\pi^2)^{-1}|z|^{-2}$ and using again
$(\mathfrak{H}_{k}^{3})$ with
$1
\leq k
\leq p-1$, we get for $\varepsilon$ sufficiently small,
\begin{eqnarray*} I^\varepsilon_{p,j}(A_{5}^{(2)},y) \leq
4\,\varepsilon^{4p}
\int_{2f(\varepsilon)}^{\infty}
\frac{dr}{r^{3}} \left(  \log \frac{r}{\varepsilon} 
\right)^{p-2},
\end{eqnarray*}
and as before in the estimate for $I^\varepsilon_{p,j}(A_3^{(1)},y)$,
it follows that   
$$ \lim_{\varepsilon \to 0}\left(
\sup_{4\varepsilon < |y| < f(\varepsilon)}
\varepsilon^{-4p} |y|^{2}\left( \log 
\frac{f(\varepsilon)}{\varepsilon} 
   \right)^{1-p} I^\varepsilon_{p,j}(A_{5}^{(2)},y) \right)=0.$$
We get $(\mathfrak{H}_{p}^{2})$ by combining the preceding estimates 
on $I^\varepsilon_{p,j}(A_i^{(2)},y)$ for $1\leq i\leq 5$, and using
(\ref{decoupyinter}).
\vskip0.2cm
\PA We now prove $(\mathfrak{H}_{p}^{3})$, and we thus assume that 
$|y|\geq f(\varepsilon)$. Let $C'>1$ be such that  $1-(C')^{-1}=(1+\beta)^{-1}$.
For $\varepsilon>0$ sufficiently small, we can 
split the integral in $I^\varepsilon_{p,j}(y)$ into five parts
corresponding to the sets
\begin{itemize}
\item
$A_{1}^{(3)}=B_{ 4\varepsilon}$\;,
\item
   $A_{2}^{(3)}= B_{\hat{f}(\varepsilon)} \setminus B_{4\varepsilon}$\;,
\item
$A_{3}^{(3)}= B_{|y|/C'} \setminus B_{\hat{f}(\varepsilon)}$\;,
\item
$A_{4}^{(3)} = B_{2|y|} \setminus B_{|y|/C'}$\;,
\item
$A_{5}^{(3)}= \mathbb{R}^{4} \setminus B_{ 2|y|}$\;.
\end{itemize}
We have then
\begin{equation} \label{decoupygrand}
\mathbb{N}_{y} \left[ \langle \z,\varphi_\varepsilon\rangle^{p} \right] = p! \sum_{j=1}^{p-1} 
\sum_{i=1}^{5} I^\varepsilon_{p,j}(A_{i}^{(3)},y).
\end{equation}
Using $(\mathfrak{H}_{k}^{1})$ with $1 \leq k \leq p-1$, we get for 
$\varepsilon$ small that
   \begin{eqnarray*}
I^\varepsilon_{p,j}(A_{1}^{(3)},y) \leq \varepsilon^{4p-4}
\left( \log  \frac{f(\varepsilon)}{\varepsilon} \right)^{p-2}
\frac{4}{|y|^{2}} \int_{0}^{4\varepsilon}r^{3}dr,
\end{eqnarray*}
so that
\begin{eqnarray*}
   \lim_{\varepsilon \to 0}\left( \sup_{|y| \geq f(\varepsilon)}
\varepsilon^{-4p} |y|^{2}
\left( \log \frac{f(\varepsilon)}{\varepsilon} \right)^{1-p}
I^\varepsilon_{p,j}(A_{1}^{(3)},y)\right) =0.
\end{eqnarray*}
Then, using
$(\mathfrak{H}_{k}^{2})$ with $1 \leq k \leq p-1$, we obtain for
$\varepsilon$ small enough,
\begin{eqnarray*}
I^\varepsilon_{p,j}(A_{2}^{(3)},y) \leq \varepsilon^{4p}
\left( \log \frac{f(\varepsilon)}{\varepsilon} \right)^{p-2}
\frac{4}{|y|^{2}} \int_{4\varepsilon}^{\hat{f}(\varepsilon)}r^{-1}dr.
\end{eqnarray*}
Thus, using (\ref{limitffhat}), we have
\begin{eqnarray*}
   \lim_{\varepsilon \to 0}\left( \sup_{|y| \geq f(\varepsilon)}
\varepsilon^{-4p} |y|^{2}
\left( \log  \frac{f(\varepsilon)}{\varepsilon}  \right)^{1-p}
I^\varepsilon_{p,j}(A_{2}^{(3)},y) \right)=0.
\end{eqnarray*}
If $z \in A_{3}^{(3)}$,
we have
$G(y,z) \leq (1-(C')^{-1})^{-2} G(0,y) =
(1+\beta)^{2}G(0,y)$.
Since $\hat{f} \in \mathcal{F}$,
we can use
$(\mathfrak{H}_{k}^{3})$ with $1 \leq k \leq p-1$ to get for
$\varepsilon$ sufficiently small,
\begin{eqnarray*}
I^\varepsilon_{p,j}(A_{3}^{(3)},y) &\leq&
(1+\beta)^2\left(\frac{\overline\varphi+\beta}{2\pi^2}\right)^p \varepsilon^{4p}
\frac{1}{|y|^{2}} \int_{\hat{f}(\varepsilon)}^{|y|/C'} \frac{dr}{r} \left(
\log \frac{r}{\varepsilon}  \right)^{p-2} \\
&=& \frac{(1+\beta)^2}{p-1}\left(\frac{\overline\varphi+\beta}{2\pi^2}\right)^p \varepsilon^{4p}
\frac{1}{|y|^{2}} \left(
\log \frac{|y|}{C' \hat f(\varepsilon)} \right)^{p-1},
\end{eqnarray*}
Recalling (\ref{limitffhat}),
we obtain that
\begin{eqnarray*}
   &&\limsup_{\varepsilon \to 0} \left(\sup_{|y| \geq f(\varepsilon)}
\varepsilon^{-4p} |y|^{2}
\left( \log  \frac{|y|}{\varepsilon} \right)^{1-p}
I^\varepsilon_{p,j}(A_{3}^{(3)},y)\right)\\
&&\qquad \leq \frac{(1+\beta)^2}{p-1}\left(\frac{\overline\varphi+\beta}{2\pi^2}\right)^p
< \frac{1}{p-1}\; \left(\frac{\overline\varphi+\beta_0}{2\pi^2}\right)^p,
\end{eqnarray*}
using our choice of $\beta$.
If $z \in A_{4}^{(3)}$, we have
$|z|^{-2} \leq {C'}^{2} |y|^{-2}$. Since
   $\hat{f} \in \mathcal{F}$, we obtain from
$(\mathfrak{H}_{k}^{3})$ with $1 \leq k \leq p-1$
that for sufficiently small $\varepsilon$,
\begin{eqnarray*}
I^\varepsilon_{p,j}(A_4^{(3)},y)
&\leq &\varepsilon^{4p}\,
\frac{{C'}^{4}}{|y|^{4}} \int_{\{|y|/C'\leq |z|\leq 2|y|\}} dz G(y,z)  \left(
\log \frac{|z|}{\varepsilon}  \right)^{p-2} \\
   &\leq&  \varepsilon^{4p}\,
\frac{{C'}^{4}}{|y|^{4}}
   \int_{\{|z|\leq 3|y|\}} \frac{dz}{2\pi^{2}|z|^{2}}  \left(
\log_+\left( \frac{|z+y|}{\varepsilon} \right) \right)^{p-2} \\
   &\leq&  \varepsilon^{4p}\,
\frac{9{C'}^{4}}{2|y|^{2}}
\left(\log \frac{4|y|}{\varepsilon}  \right)^{p-2}.
\end{eqnarray*}
It follows that
\begin{eqnarray*}
   \lim_{\varepsilon \to 0} \left(\sup_{|y| \geq f(\varepsilon)}
\varepsilon^{-4p} |y|^{2}
\left( \log  \frac{|y|}{\varepsilon} \right)^{1-p}
I^\varepsilon_{p,j}(A_{4}^{(3)},y) \right)=0.
\end{eqnarray*}
Finally, using $(\mathfrak{H}_{k}^{3})$ with $1 \leq k \leq p-1$, we get for $\varepsilon$ small,
\begin{eqnarray*}
I_{p,j}^{\varepsilon}(A_5^{(3)},y) &\leq&
4\,\varepsilon^{4p}
   \int_{2|y|}^{\infty} \frac{dr}{r^{3}} \left(
\log \frac{r}{\varepsilon} \right)^{p-2} \\
   &\leq& K\, \varepsilon^{4p}
   \frac{1}{|y|^{2}} \left(
\log \frac{2|y|}{\varepsilon} \right)^{p-2},
\end{eqnarray*}
for some constant $K$ depending only on $p$. Thus,
   \begin{eqnarray*}
   \lim_{\varepsilon \to 0} \left(\sup_{|y| \geq f(\varepsilon)}
\varepsilon^{-4p} |y|^{2}
\left( \log  \frac{|y|}{\varepsilon}  \right)^{1-p}
I^\varepsilon_{p,j}(A_{5}^{(3)},y) \right)=0.
\end{eqnarray*}
Combining our estimates on $I_{p,j}^{\varepsilon}(A_i^{(3)},y)$
for $1\leq i\leq 5$ and then summing over $j$ using
(\ref{decoupygrand}), we get that $(\mathfrak{H}_{p}^{3})$ holds
for the given $f$ and $\beta_0$.
This completes the proof of the bounds $(\mathfrak{H}_{p}^{1})$,
$(\mathfrak{H}_{p}^{2})$ and $(\mathfrak{H}_{p}^{3})$,
for every $p\geq 1$.
\vskip0.2cm
\PA It now follows from (\ref{lowerestimate}) and $(\mathfrak{H}_{p}^{3})$
that, for every $p\geq 1$,  
$$ \lim_{\varepsilon \to 0}
\varepsilon^{-4p}
\left( \log  \frac{|x|}{\varepsilon} \right)^{1-p}
\mathbb{N}_{x}\left[\langle \z,\varphi_\varepsilon\rangle^{p}\right] =  
p!\,\left(\frac{\overline\varphi}{2\pi^2}\right)^p\, |x|^{-2} ,$$
uniformly when $x$ varies over a compact subset of $\R^4\setminus\{0\}$.
This completes the proof of Lemma \ref{estimomentp}. $\qquad \Box$

\end{document}